\documentclass{article}[12pt]
\usepackage{amsmath,amsfonts,amsthm}
\usepackage{amssymb}

\usepackage{epsfig}

\usepackage{graphics}
\usepackage{graphicx}

\usepackage{latexsym}
\usepackage{graphics}
\usepackage{amsmath}
\usepackage{amsthm}
\usepackage{xspace}
\usepackage{amssymb}
\usepackage{color}
\usepackage{cite}
\usepackage{latexsym}
\usepackage{graphics}
\usepackage{amsmath}
\usepackage{amsthm}
\usepackage{xspace}
\usepackage{amssymb}
\usepackage{epsfig}

\newcommand{\beql}[1]{\begin{equation}\label{#1}}
\newcommand{\eeql}{\end{equation}}
\newcommand{\eqn}[1]{(\ref{#1})}

\newcommand{\R}{\mathbb{R}}
\newcommand{\pr}{\mathbb{P}}
\newcommand{\E}{\mathbb{E}}

\newcommand{\cx}{{\cal X}}

\newcommand{\cn}{{\cal N}}
\newcommand{\cj}{{\cal J}}

\newcommand{\Z}{\mathbb{Z}}
\newcommand{\barZ}{\bar{\Z}}

\newtheorem{thm}{Theorem}
\newtheorem{lem}[thm]{Lemma}
\newtheorem{prop}[thm]{Proposition}

\newtheorem{definition}[thm]{Definition}

\newtheorem{rem}{Remark}

\begin{document}

\title{Pull-based load distribution
\\ in large-scale heterogeneous service systems
}

\author
{
Alexander L. Stolyar \\
Lehigh University\\
200 West Packer Avenue, Room 484\\
Bethlehem, PA 18015 \\
\texttt{stolyar@lehigh.edu}
}

\date{\today}

\maketitle

\begin{abstract}

The model is motivated by the problem of load distribution
in large-scale cloud-based data processing systems.
We consider a heterogeneous
service system, consisting of multiple large server pools.
The pools are different in that their servers may have different processing speed
and/or different buffer sizes (which may be finite or infinite).
We study an asymptotic regime in which the customer arrival rate
and pool sizes scale to infinity simultaneously, in proportion 
to some scaling parameter $n$.

Arriving customers are assigned to the servers by a ``router'',
according to a {\em pull-based} algorithm, called PULL. 
Under the algorithm, each server sends a ``pull-message''
to the router, when it becomes idle; the router assigns an arriving customer
to a server
according to a randomly chosen available pull-message, if there are any,
or to a random server, otherwise.

Assuming sub-critical system load,
we prove asymptotic optimality of PULL. Namely,
as system scale $n\to\infty$, 
the steady-state probability of an arriving customer
experiencing blocking or waiting, vanishes.
We also describe some generalizations of the model and PULL algorithm,
for which the asymptotic optimality still holds.

\end{abstract}




\section{Introduction}

Modern cloud-based data processing systems are characterized 
by very large 
scale \cite{G11}. 
Service requests in such systems are processed by large-scale pools
of ``servers'', which may be physical or virtual. 
The design of efficient load distribution, i.e. routing of
arriving requests to the servers, 
in such large-scale systems poses significant challenges;
especially in {\em heterogeneous} systems, where the servers 
may have different capabilities.
Key objectives of a load distribution (routing) scheme are:
(a) keep the request response times and blocking probabilities
small and (b) keep the router/server-signaling overhead at a manageable level.

In this paper we consider a generic heterogeneous service system, 
consisting of multiple server pools.
The pools are different in that their servers may have different processing speed
and/or different buffer sizes (which may be finite or infinite).
We propose and study a {\em pull-based} routing (load distribution) algorithm, 
refered to as PULL.

The basic model and basic PULL algorithm are as follows.
(The model and the algorithm allow multiple generalizations;
see the end of this section and  
Section~\ref{sec-generalization}.)
Each new customer (service request) first arrives 
at a single ``router'' (or ``dispatcher''), which immediately
sends it to one of the
servers, as described below. 
Each server processes customers 
in the first-come-first-serve (FCFS) order.
At time instants when server becomes idle it sends a ``pull-message'' 
to the router. 
(In a more general version of PULL in Section~\ref{sec-parallel-service},
a pull-message is sent at time instants when server idleness {\em increases.})
Upon arrival of a new customer in the system,
if router has available pull-messages,
it sends the customer to one of the servers
according to an available pull-message, chosen randomly
uniformly, and ``destroys'' this pull-message.
If no pull-messages are available, the router sends 
the customer randomly uniformly to one of the servers in the 
system. We assume that pull-messages 
are never lost or ``disappear'' for any reason.
This effectively means that at any time the router ``knows'' which servers
are idle (but has {\em no other information} about the servers' parameters
or state).

We consider an asymptotic regime in which the customer arrival rate
and pool sizes scale to infinity simultaneously in proportion 
to scaling parameter $n$; we choose $n$ to be the total number of servers (in all pools)
in the system. Specifically, the arrival rate is $\lambda n$ and the server
pool sizes are $\beta_1 n, \ldots, \beta_J n$, for some positive constants
$\lambda$ and $\beta_j, ~j=1,\ldots, J$, $\sum_j \beta_j =1$;
the service rate at one server in pool $j$ is $\mu_j>0$.
We assume sub-criticality of the system load: 
$\lambda < \sum_j \beta_j \mu_j$. 

{\bf Our main result:} 
{\em PULL algorithm is asymptotically optimal; namely,
as $n\to\infty$, the steady-state probability of an arriving customer
experiencing blocking or waiting, vanishes.}

A pull-based approach to load distribution has
 been relatively recently introduced in the literature
\cite{BB08,G11}. However, a rigorous analytic study of 
pull-based algorithms is lacking,
to the best of our knowledge. 
(For example, the analysis in \cite{G11} does {\em not} imply the asymptotic optimality of PULL,
even in homogeneous systems.)
Moreover, there are no analytic studies
of pull-based schemes in heterogeneous systems, 
again, 
to the best of our knowledge.

Pull-based algorithms are very attractive for practical implementation.
Their advantages are best illustrated (see also \cite{G11}) in comparison
with the celebrated {\em power-of-d-choices},
or {\em join-shortest-queue(d)} [JSQ(d)]
 algorithm \cite{VDK96,Mitz2001,BLP2012-jsq-asymp-indep,BLP2013-jsq-asymp-tail}.
The  JSQ(d) algorithm routes an arriving customer to the server
that has the shortest queue out of the $d$ servers picked uniformly
at random. (Integer $d\ge 1$ is the algorithm parameter.)

Consider first a {\em homogeneous} system with all $n$ servers having same service rate
$\mu>0$, exponentially distributed service times, and infinite buffer sizes. 
(This is the setting 
in \cite{VDK96}.) The subcriticality condition is $\lambda<\mu$.
Denote by $p^n_k$ the steady-state probability that,
in the $n$-th system, a given server has queue length
at least $k\ge 0$. The main result of \cite{VDK96} is
\beql{eq-VDK-tail}
\lim_{n\to\infty} p^n_k = (\lambda/\mu)^{(d^k-1)/(d-1)}, ~~k\ge 0.
\end{equation}
In the case $d=1$ (which is equivalent to
random uniform routing) the RHS above is $(\lambda/\mu)^k$.
Therefore, if $d\ge 2$, the steady-state queue length tail probability
decays dramatically faster than in the case of random uniform routing.
Note that JSQ(d) does not need to maintain information on the queue lengths
at all servers. The required message exchange rate between router and the servers
is $2d$ messages per one customer. ($d$ queue length request messages
from router to servers and $d$ responses.)
To summarize, the key advantage of JSQ(d) with small $d>1$, say JSQ(2),
over random routing (JSQ(1)) is that a dramatic reduction in queue length and
waiting time is achieved at the cost of only a small message-exchange rate of
$2d$ per customer.

PULL algorithm provides further substantial improvements over JSQ(d).
Indeed, our results show that, under PULL,
$$
\lim_{n\to\infty} p^n_k = 0, ~~k\ge 2,
$$
and in fact the steady-state probability of an arriving customer having to wait
for service vanishes as well. 
The message-exchange rate of PULL in steady-state
is one message per customer.
(So, for example, this is $4$ times less that under JSQ(2).)
Therefore, when system scale $n$ is large, PULL both 
dramatically reduces (in the limit -- eliminates) queueing delays
and very substantially reduces the message-exchange rate,
compared to JSQ(d).

Suppose now that the servers have finite buffer sizes $B\ge 1$. 
For any fixed $B$, no matter how large, under JSQ(d), the steady-state blocking 
probability does {\em not} vanish as $n\to\infty$.
In contrast, under PULL, both the blocking and waiting probabilities vanish.
This is true even when $B = 1$, i.e. in the pure blocking system, where
each arriving customer either immediately goes to service or is blocked.

Further, consider heterogeneous systems, 
which are the focus of this paper. {\em In heterogeneous systems 
the JSQ(d) algorithm is not even appropriate in general.} 
To illustrate, suppose there are two server pools, each of size $n/2$,
with service rate parameters $\mu_1=2$ and $\mu_2=1/3$.
Assume infinite buffer sizes at all servers.
The arrival rate is $n$, so that the subcriticality holds: $1 < (1/2)2+(1/2)(1/3)$.
Under JSQ(2) this system is {\em unstable}, because the second (slower) 
pool will receive new arrivals at the rate at least $(1/4)n$,
while its maximum service rate is $(1/2)n(1/3)$.
In contrast, under PULL, the system is stable (for sufficiently large $n$)
and the probability of waiting vanishes, as our results show.

Finally, we remark that our basic model and the PULL algorithm
can be easily generalized, so that the asymptotic optimality 
of (more general) PULL still holds -- essentially same proofs
as for the basic model work. In Section~\ref{sec-generalization}
 we discuss two 
such generalizations: (a) for the case when a server processing
rate depends on the queue length and (b) for more general
service time distributions, namely, those with decreasing hazard rate (DHR).

\subsection{Brief literature review and summary of contributions} 

The literature on load distribution in service systems is extensive; see e.g.
\cite{BLP2013-jsq-asymp-tail,G11} for 
good up-to-date overviews.
A lot of previous work is focused on {\em load balancing}, which, we note, 
is only one of  possible objectives of {\em load distribution}.
The PULL algorithm, studied in this paper, does not attempt and does not in general achieve load balancing
in the sense of equal load of the servers. (It does provide load balancing {\em within each server pool}.)
Nevertheless, it achieves the asymptotic optimality in the sense of eliminating customer waiting and blocking.

The JSQ(d) algorithm, for homogeneous systems, has 
received much attention, since it was introduced in the seminal work \cite{VDK96}.
(See \cite{Mitz2001,BLP2012-jsq-asymp-indep,BLP2013-jsq-asymp-tail}
for reviews.) Paper \cite{VDK96} considers the (homogeneous) system
with exponential service time distribution, under the same asymptotic regime as in this paper,
and proves the limit \eqn{eq-VDK-tail}
for queue length distribution. Significant generalizations of the results of \cite{VDK96}
are obtained in \cite{BLP2012-jsq-asymp-indep,BLP2013-jsq-asymp-tail};
in particular, these papers establish the queue length distribution limit
for the case when the service time distribution has decreasing hazard rate (DHR).

The basic idea of a pull-based load distribution is to make servers ``pull'' 
customers for service, as opposed to router ``pushing'' it to them (as in JSQ(d)).
Paper \cite{BB08} proposes various pull-based schemes, with the focus on practical use,
and studies them via simulation. Recent work \cite{G11} considers 
a pull-based algorithm in a homogeneous system;
the model in  \cite{G11} is more general than ours in that it has {\em multiple} routers,
each handling equal fraction of customer arrivals;
the analytic and simulation study in the paper shows potentially
significant advantages of a pull-based approach over JSQ(d).
(But, as we mentioned, it does not prove asymptotic optimality.)

{\bf Summary of this paper contributions:} \\(1) We propose a 
specific pull-based load distribution algorithm, called PULL.
\\(2) We rigorously prove the asymptotic optimality of PULL 
(namely, elimination of waiting and blocking)
in a heterogeneous service system. In particular, this proves 
that PULL asymptotic performance is much superior to that of the celebrated
power-of-d-choices [JSQ(d)] algorithm.
\\(3) We present two generalizations of the model and the PULL algorithm,
for which asymptotic optimality prevails: for the queue length dependent service rates
and for service time distributions with DHR.

\subsection{Basic notation}

Symbols $\R, \R_+, \Z, \Z_+$ denote the sets of real, real non-negative,
integer, and integer non-negative numbers, respectively.
For finite- or infinite-dimensional vectors, the vector inequalities 
are understood component-wise. We write simply $0$ for a zero-vector.
We use notation $x(\cdot)=(x(t), ~t\ge 0)$ for both a random process
and its realizations, the meaning is determined by the context;
the state space (of a process) and the metric and/or 
topology on it are defined where appropriate,
and we always consider Borel $\sigma$-algebra on the state space.
Abbreviation {\em u.o.c.} means 
{\em uniform on compact sets} convergence, and {\em w.p.1} means
{\em with probability 1.} Notations $\Rightarrow$ and $\stackrel{d}{=}$
signify convergence and equality {\em in distribution}, respectively,
for random elements.
For a process $x(\cdot)$, we denote by $x(\infty)$ a random element
whose distribution is {\em the lower invariant measure} of the process
(defined formally in the text); if the process has unique
stationary distribution, it is equal to the lower invariant measure.
For $a\in \R$, $\lfloor a \rfloor$ denotes the largest integer less than 
or equal to $a$.

\subsection{Layout of the rest of the paper} 

The formal model, asymptotic regime, PULL algorithm definition
and the main result (Theorem~\ref{thm-infinite})
are given in Section~\ref{sec-model}.
In Section~\ref{sec-monotonicity} we study
properties of the underlying Markov process,
related to -- and stemming from -- its monotonicity.
Fluid limits (as $n\to\infty$) of the process
are studied in Section~\ref{sec-fluid-many-server}.
The proof of Theorem~\ref{thm-infinite} is given
in Section~\ref{sec-proof-infinite}.
In Section~\ref{sec-generalization} we discuss 
generalizations of the model and PULL algorithm,
for which our main results still hold, with essentially same proofs.

\section{Model and main result}
\label{sec-model}

\subsection{Model structure}
\label{subsec-model}

Customers for service arrive according to a Poisson 
process of rate $\Lambda>0$.
There are $J\ge 1$ server pools. Pool $j \in \cj \equiv \{1,\ldots,J\}$
consists of $N_j$ identical servers.
Servers in pool 1 at indexed by $i \in \cn_1 =\{1,\ldots, N_1\}$,
in pool 2 by $i \in \cn_2 =\{N_1+1,\ldots, N_1+N_2\}$, and so on;
$\cn = \cup \cn_j$ is the set of all servers.
Each arriving customer is immediately routed for service to one of the servers;
the service time of a customer at a server in pool $j$
is an independent, exponentially distributed random variable 
 with mean $1/\mu_j\in (0,\infty)$, $j \in \cj$.
We assume that the customers at any server are served in the 
first-come-first-serve (FCFS) order. 
(That is, at any time
only the head-of-the-line customer at each server is served.)
The buffer size (maximum queue length) at any server in pool $j$ is $B_j\ge 1$;
we allow the buffer size to be either finite, $B_j < \infty$, or infinite, $B_j=\infty$.
A new customer, routed to a server $i\in \cn_j$, joins the queue at that server,
unless $B_j$ is finite and the queue length $Q_i=B_j$ -- in this case 
the customer is lost 
(i.e., leaves the system immediately, without receiving any service).

\begin{rem}
\label{rem-generalizations}
In the model described above, the FCFS assumption is not important
{\em as far the queue lengths in the system are concerned} -- 
any non-idling work-conserving discipline will produce the same queue
length process. In Section~\ref{sec-generalization}
 we will discuss several generalizations 
of the above model, for which our main results still hold.
Some of these generalizations, specifically those involving 
more general service time distributions (Section~\ref{sec-dhr}), 
{\em do} require the FCFS assumption.
\end{rem}

\subsection{Asymptotic regime}
\label{sec-asymptotic-reg}

We consider the following (many-servers) asymptotic regime. 
The total number of servers $n=\sum_j N_j$ is the scaling parameter, which
increases to infinity; the arrival rate and the server pool sizes increase in proportion 
to $n$, $\Lambda= \lambda n$, $N_j=\beta_j n, ~j\in \cj$, where 
$\lambda, \beta_j, j\in \cj$, are positive constants, $\sum_j \beta_j = 1$.
(To be precise, the values of $N_j$ need to be integer, 
e.g. $N_j=\lfloor \beta_j n \rfloor$. Such definition
would not cause any problems, besides clogging
notation, so we will simply assume that all 
$\beta_j n$ ``happen to be'' integer.)
We assume that the subcritical load condition holds:
\beql{eq-load}
\lambda < \sum_j \beta_j \mu_j.
\end{equation}

\subsection{PULL routing algorithm}
\label{sec-alg}

We study the following pull-based algorithm.
\begin{definition}[PULL algorithm]
\label{def-pull-basic}
At any given time the algorithm (router) has exactly
one {\em pull-message} from each idle server (i.e., server with zero queue length) in the system.
(In other words, the algorithm ``knows'' which servers are idle.)
Each arriving customer is routed immediately to one of the servers. 
If there are available pull-messages (idle servers), the customer is
routed to one of the idle servers, chosen randomly uniformly.
If there are no available pull-messages (idle servers), the customer is 
routed to one of the servers in the system, chosen randomly uniformly.
\end{definition}

A practical implementation of PULL algorithm (which motivates it name) is as follows.
Assume that pull-messages are never lost.
When a server is ``initialized'', it sends one pull-message to the router.
After that, the server sends one new pull-message to the router immediately after
any service completion that leaves the server idle.
When a customer arrives, the router picks one of the available pull-messages
uniformly at random,  sends the customer to the corresponding server,
and destroys the pull-message.
If router has no available pull-messages when a customer arrives, it sends the customer to
one of the servers, chosen uniformly at random.
Thus, the algorithm is easily implementable.
Of course, in the algorithm {\em analysis}, there is no need to consider 
the pull-message mechanism -- we just assume that the current set of idle servers is known at any time.

We will discuss implementation aspects of PULL in more detail in Section~\ref{sec-practical}, after formally stating our main result.

\subsection{Main result}
\label{sec-result}


In the system with parameter $n$, the system state is the vector
$Q^n=(Q^n_i, ~i\in \cn)$, where $Q^n_i \in \Z_+$ is the
queue length at server $i$. 

Due to symmetry of servers within each pool,
the alternative -- {\em mean field}, or {\em fluid-scale}  
-- representation of the process is
as follows. 
Define $x^n_{k,j}$ as the fraction of the (total number of) servers,
which are in pool $j$ and have queue length greater than or equal to $k$. 
We consider 
$$
x^n = (x^n_{k,j}, ~k\in \Z_+, ~j\in \cj),
$$
to be the system state, and will view states
$x^n$, for any $n$, as elements of the common space
$$
\cx = \{x = (x_{kj}, ~k\in \Z_+, ~j\in \cj)~|~ \beta_j = x_{0j} \ge x_{1j} \ge x_{2j} \ge \cdots
\ge 0\},
$$
equipped with metric
\beql{eq-metric}
\rho(x,x') = \sum_{j} \sum_k 2^{-k} \frac{|x_{kj}-x'_{kj}|}{1+|x_{kj}-x'_{kj}|},
\end{equation}
and the corresponding
Borel $\sigma$-algebra. Space $\cx$ is compact. 

For any $n$, the process $Q^n(t), ~t\ge 0,$ -- and its projection $x^n(t), ~t\ge 0,$ -- 
is a continuous-time, 
countable state space, irreducible Markov process. 
(For any $n$, the state space of $x^n(\cdot)$ is a countable subset of $\cx$.)
If the buffer sizes $B_j$ are finite in {\em all} pools $j$, 
the state space is obviously finite, and therefore the process $Q^n(\cdot)$
(and then $x^n(\cdot)$)
is ergodic, with unique stationary distribution. We will prove 
(in Theorem~\ref{thm-infinite}) that, in fact, the ergodicity 
holds in the general case, when buffer sizes $B_j$ may be infinite in some 
or all pools.

Define numbers $\nu_j \in (0,\beta_j)$, $j\in \cj$,
uniquely determined by the conditions
\beql{eq-invar-point1}
\lambda = \sum_j \nu_j \mu_j, 
~~~ \nu_j \mu_j / (\beta_j - \nu_j) = \nu_\ell \mu_\ell / (\beta_\ell - \nu_\ell),
~\forall j,\ell\in\cj.
\end{equation}
Let us define the {\em equilibrium point} $x^* \in \cx$ by
\beql{eq-invar-point2}
x^*_{1,j} = \nu_j, ~~ x^*_{k,j} = 0,~k \ge 2, ~~~ j\in\cj.
\end{equation}

The meaning of the equilibrium point $x^*$ definition in \eqn{eq-invar-point1}-\eqn{eq-invar-point2}
is clear. Point $x^*$ is such that the fraction $\nu_j < \beta_j$ of servers
(out of the total number of servers) in pool $j$ is occupied by exactly one customer, while the remaining servers in pool $j$ are idle.
The numbers $\nu_j$ are (uniquely) determined by the condition \eqn{eq-invar-point1}, which simply says that the rate at which 
new arrivals are routed to pool $j$ (it is proportional to $\beta_j - \nu_j$) is equal to the service-completion/departure rate from pool $j$ 
(it is proportional to $\nu_j \mu_j$).

If the process $x^n(\cdot)$ is ergodic, it has unique stationary distribution;
in this case, we denote by $x^n(\infty)$ a random element with the distribution
equal to the stationary distribution of $x^n(\cdot)$.
(In other words, $x^n(\infty)$ is a random process state in stationary regime.)
Our main result is the following 
\begin{thm}
\label{thm-infinite}
For all sufficiently large $n$, the Markov process $x^n(\cdot)$ is ergodic
(and then has a unique stationary distribution), and $x^n(\infty)\Rightarrow x^*$.
\end{thm}

Given the definition of $x^*$, the result implies that, as $n\to\infty$, 
the steady-state probability of having an idle server in the system,
goes to $1$. Consequently, 
the steady-state probability of an arriving customer
experiencing blocking or waiting, vanishes.

\subsection{Discussion of implementation aspects of PULL algorithm}
\label{sec-practical}

\subsubsection{The notion of servers pools is purely logical.}

Note that PULL algorithm uses only the information about which servers are idle;
it needs to know {\em neither} the queue lengths at the servers 
(besides it being zero or not), {\em nor}
their processing speed (i.e. which pool $j$ they belong to), {\em nor}
their buffer sizes. In other words, from the ``point of view'' 
of the router, all servers form a single pool, and the router
need {\em not} know anything about the servers, besides them being currently 
idle or not.

This in particular means that our model's notion of server pools, each consisting of identical servers of a certain type, is {\em purely logical},
used for the purposes of analysis only. A real system may consist of a single or multiple pools of {\em non-identical} 
servers. In this case, we consider all servers of a particular type as forming a {\em logical} pool.
Our results still apply, as long as the number of servers of each type in the entire system is large.

\subsubsection{Pull-message mechanism.} 

We already mentioned that PULL algorithm very substantially reduces the message exchange between the router and the servers. (It is $2d$ times less than that of JSQ(d).) Furthermore, pull-messages do {\em not} contribute
to the routing delay: an arriving customer in not waiting at the router for any pull-message, the routing decision is made 
immediately, based on the pull-messages
currently available. This is unlike the JSQ(d) algorithm, where each arriving customer waits for the queue-length request/response message exchange to complete, before being routed. (See also \cite{G11} for a discussion of this issue.)

It may appear that a disadvantage of PULL algorithm, compared to JSQ(d), is that the router needs to maintain the list
of available pull messages. In fact, this issue is insignificant for the following reason. Under {\em any} routing algorithm,
including JSQ(d) and PULL, the router needs to 
have the list of all servers in the system. (It needs to know server ``addresses", in whatever form, to do actual routing of customers.)
A pull-message availability (or not) for a given server, adds just one bit to the server's entry on the list.
For the purposes of making the random choice of an available pull-message efficient, it might be beneficial to connect the corresponding
server entries to form a virtual list; even in this case, this just adds one additional field to each server entry.

\subsubsection{Amount of computation.}

The JSQ(d) algorithm needs to generate $d$ random (or pseudo-random) numbers per each routed customer.
Under PULL algorithm, only one random number is generated per each customer.

\section{More general view of the process. 
Monotonicity. Lower invariant measure}
\label{sec-monotonicity}

All results in this section concern a system with a fixed $n$.

It will be convenient to consider a more general system and the Markov process.
Namely, we assume that the queue length in any server $i\in\cn_j$
within a pool $j$ with infinite buffer size ($B_j=\infty$), can be 
infinite. In other words, $Q_i(t)$ can take values in the
set $\barZ_+ \doteq \Z_+ \cup \{\infty\}$, which is the
one-point compactification of $\Z_+$, containing the ``point at infinity.''
We consider the natural topology and order relation on $\barZ_+$.
Obviously, $\barZ_+$ is compact. (Note that if $A$ is a finite subset
of $\Z_+$, then sets $A$ and $\barZ_+ \setminus A$ are both closed and open.)

Therefore, the state space of the generalized version of 
Markov process $Q^n(\cdot)$ is the compact set $\barZ_+^n$. 
The process transitions are defined in
exactly the same way as before, with the additional convention that
if $Q^n_i(t)=\infty$, then neither new arrivals into this queue nor
service completions in it, change the infinite queue length value,
and therefore $Q^n_i(\tau)\equiv \infty$ for all $\tau\ge t$.

The corresponding generalized version of 
the process $x^n(\cdot)$ is defined as before; if at time $t$ some of the
queues in pool $j$ are infinite, then $x^n(t)$ is such that
$\lim_{k\to\infty} x^n_{k,j}(t) >0$. Note that the 
state space of the generalized $x^n(\cdot)$ is still the
compact set $\cx$, as defined above.

It is easy to see that, for each $n$, the (generalized versions of) 
processes $Q^n(\cdot)$ and $x^n(\cdot)$ are Feller continuous.

Vector inequalities, $Q'\le Q''$ for $Q', Q'' \in \barZ_+^n$
and $x'\le x''$ for $x', x'' \in \cx$,
are understood component-wise.
The stochastic order relation $Q'\le_{st} Q''$ [resp. $x' \le_{st}  x''$] for
{\em random elements}  taking values in $\barZ_+^n$ [resp. $\cx$]
means that they can be constructed on the
same probability space so that $Q'\le Q''$ [resp. $x' \le  x''$] holds w.p.1.

For any $n$, the processes $Q^n(\cdot)$ and $x^n(\cdot)$ are {\em  monotone}.
Namely, the following property holds.
(For a general notion of monotonicity cf. \cite{Liggett-book}.)

\begin{lem}
\label{lem-monotone}
Consider two version of the process, $Q^n(\cdot)$ and $\bar Q^n(\cdot)$ 
[resp. $x^n(\cdot)$ and $\bar x^n(\cdot)$],
with fixed initial states $Q^n(0) \le \bar Q^n(0)$ [resp. $x^n(0) \le \bar x^n(0)$].
Then, the processes can be constructed on a common probability space,
so that, w.p.1, $Q^n(t) \le \bar Q^n(t)$ [resp. $x^n(t) \le \bar x^n(t)$]
for all $t\ge 0$. Consequently, 
$Q^n(t) \le_{st} \bar Q^n(t)$ [resp. $x^n(t) \le_{st} \bar x^n(t)$]
for all $t\ge 0$.
\end{lem}

{\em Proof.} It suffices to prove the result for $Q^n(\cdot)$ and $\bar Q^n(\cdot)$.
We will refer to the systems, corresponding to $Q^n(\cdot)$ and $\bar Q^n(\cdot)$,
as ``smaller'' and ``larger'', respectively.
It is clear how to couple 
the service completions in the two systems,
so that any service completion preserves the
$Q^n(t) \le \bar Q^n(t)$ condition. We make the arrival process to be common for both systems.
It suffuces to show that condition $Q^n(t) \le \bar Q^n(t)$ is preserved 
after any arrival. Suppose the (joint) system state just before a customer arrival is such that
$Q^n \le \bar Q^n$. If all servers in both the smaller and larger system are busy, 
we make a {\em common} random uniform assignment of the arrival to one of the servers.
If all servers are busy in the larger system, but there are idle servers in the smaller one,
we make {\em independent} assignments in the two systems, according to the algorithm.
In the case when there are idle servers in both systems, obvoiusly the idle servers in the larger system
form a subset of those in the smaller one. Then we do the following.
We make uniform random choice of an idle server in the smaller system, we assign the arrival 
to that server in the smaller system, and in the larger system as well {\em as long as 
it happens to be idle in the larger system}; if that server is busy in the larger system, we do 
an additional step and assign it uniformly randomly to an idle server in the larger system.
Clearly, condition  $Q^n(t) \le \bar Q^n(t)$ is preserved in each case, and the procedure conforms 
to the PULL algorithm in both systems.
$\Box$

If the system starts from idle initial state, i.e. $Q^n(0)=0$ 
[equivalently, $x^n_{1j}(0)=0, ~j\in\cj$], then by Lemma~\ref{lem-monotone}
the process
is stochastically non-decreasing in time
\beql{eq-zero-monotone}
Q^n(t_1) \le_{st} Q^n(t_2), ~~\mbox{[resp. $x^n(t_1) \le_{st} x^n(t_2)$]}, 
~~~0\le t_1 \le t_2 < \infty.
\end{equation}
Since the state space $\barZ_+^n$ [resp. $\cx$] is compact, we 
must have convergence in distribution
$$
Q^n(t) \Rightarrow Q^n(\infty), 
~~\mbox{[resp. $x^n(t) \Rightarrow x^n(\infty)$]}, 
~~~ t\to\infty,
$$
where the distribution of $Q^n(\infty)$ [resp. $x^n(\infty)$] 
is the {\em lower invariant measure} of process 
$Q^n(\cdot)$ [resp. $x^n(\cdot)$]. (The lower invariant measure
is a stationary distribution of the process, stochastically dominated
by any other stationary distribution. Cf. \cite{Liggett-book},
in particular Proposition I.1.8(d).) 

Observe that the process $Q^n(\cdot)$ [resp. $x^n(\cdot)$], 
{\em as originally defined (without infinite queues)},
is ergodic if and only if $Q^n(\infty)$ [resp. $x^n(\infty)$] 
is proper in the sense that 
$$
\pr\{Q_i^n(\infty)<\infty,~\forall i\} = 1~~~ 
\mbox{[resp. $\pr\{x^n_{\infty,j}(\infty)=0,~\forall j\} = 1$]},
$$
where $x^n_{\infty,j}(\infty)\doteq \lim_{k\to\infty} x^n_{kj}(\infty)$.
And if the original process is ergodic, the lower invariant measure is its unique
stationary distribution. 

\begin{lem}
\label{lem-unstable-measure}
Suppose
for some $j$ and some $i\in\cn_j$, 
$$
\pr\{Q_i^n(\infty)=\infty\} >0 ~~~ 
\mbox{[and then $\pr\{x^n_{\infty,j}(\infty)>0\} > 0$]}.
$$
Then, necessarily, a stronger condition holds:
\beql{eq-instability-cond}
\pr\{Q_i^n(\infty)=\infty, ~ \forall i\in\cn_j\} =1 ~~~ 
\mbox{[and then  $\pr\{x^n_{\infty,j}(\infty)=\beta_j \} =1$]}.
\end{equation}
\end{lem}

{\em Proof.} 
Consider a stationary version of $Q^n(\cdot)$, with stationary distribution
being the lower invariant measure. Namely, $Q^n(0) \stackrel{d}{=} Q^n(\infty)$,
and then $Q^n(t) \stackrel{d}{=} Q^n(\infty)$ for all $t\ge 0$.
By the lemma assumption,
 $\pr\{Q_i^n(0)=\infty\} = \pr\{Q_i^n(\infty)=\infty\} = \delta \in (0,1]$.
By the (generalized) process definition, under the condition
$Q_i^n(0)=\infty$, w.p.1 $Q_i^n(t)\equiv \infty$ for all $t$.
By monotonicity, the process conditioned on {\em any} fixed initial state stochastically dominates 
the process starting from the idle state. Then, 
for any $k\in \Z_+$,
$$
\liminf_{t\to\infty} \pr\{Q_i^n(t) \ge k ~|~ Q_i^n(0)<\infty\} \ge 
\liminf_{t\to\infty} \pr\{Q_i^n(t) \ge k ~|~ Q^n(0)=0\} \ge
\pr\{Q_i^n(\infty) \ge k\} \ge \delta.
$$
(Recall that $\{Q_i \ge k\}$ is an open subset of $\barZ_+^n$.)
Therefore, for the overall probability (assuming $Q^n(0) \stackrel{d}{=} Q^n(\infty)$),
$$
\liminf_{t\to\infty} \pr\{Q_i^n(t) \ge k\} \ge 
\liminf_{t\to\infty} \pr\{Q_i^n(t) \ge k ~|~ Q_i^n(0)=\infty\} +
\liminf_{t\to\infty} \pr\{Q_i^n(t) \ge k ~|~ Q_i^n(0)<\infty\}
\ge \delta + (1-\delta)\delta.
$$
From here,
$$
\pr\{Q_i^n(\infty) \ge k\} \ge 
\limsup_{t\to\infty} \pr\{Q_i^n(t) \ge k\} \ge \delta + (1-\delta)\delta.
$$
(Recall that $\{Q_i \ge k\}$ is also a closed subset of $\barZ_+^n$.)
Then,
$$
\delta = \pr\{Q_i^n(\infty) =\infty \} = \lim_{k\to\infty} \pr\{Q_i^n(\infty) \ge k\}
\ge \delta + (1-\delta)\delta.
$$
We see that $(1-\delta)\delta \le 0$ and, as assumed, $\delta \in (0,1]$.
This implies $\delta=1$, that is $\pr\{Q_i^n(\infty)=\infty\}=1$. By symmetry, this is true for all servers in $\cn_j$.
$\Box$

By Lemma~\ref{lem-unstable-measure}, the non-ergodicity (instability) of the original
process is equivalent to condition \eqn{eq-instability-cond}
holding for at least one $j$.

In the rest of the paper, for a state $x^n(t)$ 
(with either finite $t\ge 0$ or $t=\infty$), we denote by
\beql{eq-x-inf}
x^n_{\infty,j}(t)\doteq \lim_{k\to\infty} x^n_{k,j}(t)
\end{equation}
the fraction of queues that are in pool $j$ and are infinite.
(Note that $x^n_{\infty,j}(t)$ is a function, but {\em not a component},
of $x^n(t)$.) Also, denote by $y^n_{k,j}(t)$ the fraction of queues that are in pool $j$ 
and have queue size {\em exactly} $k \in \barZ_+$:
\beql{eq-y-def}
y^n_{k,j}(t) \doteq x^n_{k,j}(t) - x^n_{k+1,j}(t), ~~k\in \Z_+,
\end{equation}
\beql{eq-y-def-inf}
y^n_{\infty,j}(t) \doteq x^n_{\infty,j}(t) = \lim_{k\to\infty} x^n_{k,j}(t).
\end{equation}

\section{Fluid limits}
\label{sec-fluid-many-server}

In this section, we consider limiting behavior of the sequence of process
$x^n(\cdot)$ as $n\to\infty$. (We will only consider the mean-field 
process
$x^n(\cdot)$, because this will be sufficient for proving Theorem~\ref{thm-infinite}.) 
In particular, we will define fluid sample paths (FSP), 
which arise as limits of the 
(fluid-scaled) 
trajectories $x^n(\cdot)$
as $n\to\infty$.

Without loss of generality, assume that
the Markov process $x^n(\cdot)$ for each $n$ is driven by a
common set of primitive processes, as defined next.

Let $A^n(t), ~t\ge 0$, denote the number of exogenous arrival into the system
in the interval $[0,t]$. Assume that 
\beql{eq-driving-arr}
A^n(t) = \Pi^{(a)}(\lambda n t),
\end{equation}
where $\Pi^{(a)}(\cdot)$ is an independent unit rate Poisson process.
The functional strong law of large numbers (FSLLN) holds: w.p.1
\beql{eq-flln-poisson-arr}
\frac{1}{n}\Pi^{(a)}(nt) \to t, ~u.o.c.
\end{equation}
Denote by $D^n_{k,j}(t), ~t\ge 0, ~1 \le k < \infty$, 
the total number of departures in $[0,t]$
from servers in pool $j$ with queue length $k$;
assume
\beql{eq-driving-dep}
D^n_{k,j}(t) = \Pi^{(d)}_{k,j} \left(\int_0^t n y_{k,j}^n(s) \mu_j ds\right),
\end{equation}
where $\Pi^{(d)}_{k,j}(\cdot)$ are independent unit rate Poisson processes.
(Recall that departures from -- and arrivals to -- infinite queues can be ignored,
in the sense that they do not change the system state.)
We have: w.p.1
\beql{eq-flln-poisson-dep}
\frac{1}{n}\Pi^{(d)}_{k,j}(nt) \to t, ~u.o.c., ~~~ 1 \le k < \infty.
\end{equation}

The random routing of new arrivals is constructed as follows.
There are two sequences of i.i.d. random variables, 
$$
\xi(1), \xi(2), \ldots, ~~~~\mbox{and}~~~ \zeta(1), \zeta(2), \ldots,
$$
uniformly distributed in $[0,1)$. The routing of the $m$-th arrival into the system
is determined by the values of r.v. $\xi(m)$ and $\zeta(m)$,
as follows. (We will drop index $m$, because we consider one arrival.)
Let $x^n$ denote the system state just before the arrival.
If $\sum_j y^n_{0,j} = 0$, i.e. there are no idle servers, the routing is determined
by $\zeta$ as follows.
The customer is sent to a server with $k$, $k\ge 1$, customers in pool 1,
if $\zeta \in [x^n_{k+1,1},x^n_{k,1})$, and to a server with $k=\infty$ customers in pool 1, if
$\zeta \in [0,x^n_{\infty,1})$; the customer is sent to a server with $k$, $k\ge 1$, customers in pool 2,
if $\zeta \in [\beta_1+x^n_{k+1,2},\beta_1 + x^n_{k,2})$, 
and to a server with $k=\infty$ customers in pool 2, if
$\zeta \in [\beta_1, \beta_1 + x^n_{\infty,2})$; and so on.
If $\sum_j y^n_{0,j} > 0$, i.e. there are idle servers, the routing is determined
by $\xi$ as follows.
Let $a=\sum_j y^n_{0,j}$, $p_j = y^n_{0,j}/a$.
If $\xi \in [0,p_1)$, the customer is routed 
to pool 1; if $\xi \in [p_1,p_1+p_2)$ -- to pool 2; and so on.

Denote
$$
f^n(s,u) \doteq \frac{1}{n} \sum_{m=1}^{\lfloor ns \rfloor} I\{\xi(m) \le u\}, ~~
g^n(s,u) \doteq \frac{1}{n} \sum_{m=1}^{\lfloor ns \rfloor} I\{\zeta(m) \le u\}, ~~
$$
where $s\ge 0$, $0\le u < 1$.  Obviously, from the strong law of large numbers
 and the monotonicity of $f^n(s,u)$ and $g^n(s,u)$ on both arguments, we 
have the FSLLN: w.p.1
\beql{eq-flln-random}
f^n(s,u) \to su, ~~g^n(s,u) \to su, 
~~~\mbox{u.o.c.}
\end{equation}
It is easy (and standard) to see that, for any $n$, w.p.1,
 the realization of the process $x^n(\cdot)$ 
is uniquely determined by the
initial state $x^n(0)$
and the realizations of the driving processes
$\Pi^{(a)}(\cdot)$, $\Pi^{(d)}_{k,j}(\cdot)$, $\xi(\cdot)$ and $\zeta(\cdot)$.

A set of uniformly Lipschitz continuous functions
$x(\cdot)=[x_{k,j}(\cdot),~~~ k\in \Z_+, ~j\in \cj]$
on the time interval $[0,\infty)$ we call a {\em fluid sample path} (FSP), if there exist
realizations of the primitive driving processes,
 satisfying conditions \eqn{eq-flln-poisson-arr}, \eqn{eq-flln-poisson-dep}
and \eqn{eq-flln-random}
and a fixed subsequence of $n$, along which
\beql{eq-fsp-def}
x^n(\cdot) \to x(\cdot), ~~~u.o.c.
\end{equation}

Note that, given the metric \eqn{eq-metric}
on $\cx$, condition \eqn{eq-fsp-def} is equivalent to
component-wise convergence:
$$
x^n_{k,j}(\cdot) \to x_{k,j}(\cdot), ~~~u.o.c., ~~k\in \Z_+, ~j\in \cj.
$$

For any FSP, almost all points $t\ge 0$ (w.r.t. Lebesgue measure)
are {\em regular}, namely all component functions have proper 
(equal right and left) derivatives
$(d/dt) x_{k,j}(t)$. Note that $t=0$ is {\em not} a regular point;
expression $(d/dt) x_{k,j}(t)$ for $t=0$ means right derivative (if it exists).

Analogously to notation in \eqn{eq-x-inf} - \eqn{eq-y-def-inf}, we will denote:
$$
x_{\infty,j}(t)\doteq \lim_{k\to\infty} x_{k,j}(t)
$$
$$
y_{k,j}(t) \doteq x_{k,j}(t) - x_{k+1,j}(t), ~~k\in \Z_+,
$$
$$
y_{\infty,j}(t) \doteq x_{\infty,j}(t) = \lim_{k\to\infty} x_{k,j}(t).
$$
For two FSPs $x(\cdot)$ and $\bar x(\cdot)$,
$x(\cdot) \le \bar x(\cdot)$ will mean $x(t) \le \bar x(t), ~t\ge 0$.

\begin{lem}
\label{lem-ms-fsp-conv}
Consider a sequence in $n$ of processes $x^n(\cdot)$ with deterministic initial states
$x^n(0) \to x(0)\in \cx$.
Then w.p.1 any subsequence of $n$ has a further subsequence, along which
$$
x^n(t) \to x(t) ~~~u.o.c.,
$$
where $x(\cdot)$ is an FSP.
\end{lem}

{\em Proof} is fairly standard. 
Denote by $A^n_{k,j}(t), ~k\in \Z_+, ~t\ge 0$, 
the total number of arrivals in $[0,t]$
into servers in pool $j$ with queue length $k$. 
(Recall that arrivals to infinite queues can be ignored.)
Obviously, for any $0 \le t_1 \le t_2 < \infty$
$$ 
\sum_j \sum_{1 \le k < \infty} [A^n_{k,j}(t_2)-A^n_{k,j}(t_1)] \le A^n(t_2)-A^n(t_1).
$$ 
In addition to $x^n_{k,j}(\cdot)$ (and $y^n_{k,j}(\cdot)$),
which are fluid-scaled quantities,
we define the corresponding ones for the arrival and departure processes:
$$
a^n_{k,j}(t) = \frac{1}{n} A^n_{k,j}(t),~~~ 0 \le k < \infty,
$$
$$
d^n_{k,j}(t) = \frac{1}{n} D^n_{k,j}(t),~~~ 1 \le k < \infty.
$$
All processes $a^n_{k,j}(\cdot)$ and $d^n_{k,j}(\cdot)$ are non-decreasing.
W.p.1 the primitive processes satisfy the FSLLN \eqn{eq-flln-poisson-arr},
\eqn{eq-flln-poisson-dep} and \eqn{eq-flln-random}. 
From here it is easy to observe the following:
w.p.1 any subsequence of $n$ has a further subsequence along which
the u.o.c. convergences
$$
a^n_{k,j}(\cdot) \to a_{k,j}(\cdot), ~~d^n_{k,j}(\cdot) \to d_{k,j}(\cdot),
$$
hold for all pairs $(k,j)$, where the limiting functions
$a_{k,j}(\cdot)$ and $d_{k,j}(\cdot)$ are non-decreasing,
uniformly Lipschitz continuous. The result easily follows; we omit further details.
$\Box$

\begin{lem}
\label{lem-ms-fsp-prop}
(i) If $x(\cdot)=(x(t), ~t\ge 0)$ is an FSP, then for any $\tau \ge 0$,
the time shifted trajectory $\theta_\tau x(\cdot) \doteq (x(\tau+t), ~t\ge 0)$
is also an FSP.
\\ (ii) For an FSP $x(\cdot)$,
at any $t\ge 0$, such that $\sum_j y_{0j}(t) >0$,
all derivatives $(d/dt) x_{k,j}(t)$ exist (for $t=0$, right derivatives exist)
and 
\beql{eq-deriv-1}
(d/dt) x_{1,j}(t) =  \lambda y_{0j}(t) /(\sum_{\ell}  y_{0\ell}(t)) - \mu_j  y_{1,j}(t),~ j\in \cj,
\end{equation}
\beql{eq-deriv-2}
(d/dt) x_{k,j}(t) = - \mu_j  y_{k,j}(t) \le 0, ~~ 2\le k < \infty, ~ j\in \cj.
\end{equation}
\\ (iii) If initial condition $x(0)$ of an FSP is such that $\sum_j y_{0,j}(0) >0$
and $\sum_j x_{2,j}(0)=0$,
then the FSP is unique in the interval $[0,\tau)$, where $\tau$ is the smallest time $t$
when $\sum_j y_{0,j}(t) =0$; $\tau=\infty$ if such $t$ does not exist.
\\ (iv) The FSP $x(\cdot)$ with initial condition $x(0)=x^*$ is unique, 
and it is stationary, $x(t) \equiv x^*$.
\\ (v) The FSP $x(\cdot)$ with idle initial condition, 
$x_{1,j}(0)=0, \forall j$,
is unique, monotonically increasing, $x(t_1)\le x(t_2), ~t_1 \le t_2$,
and is such that $x(t) \to x^*$. This FSP is a lower bound 
of any other FSP $\bar x(\cdot)$: $x(\cdot) \le \bar x(\cdot)$.
\\ (vi) For any $\epsilon>0$, there exist $\tau>0$ and $\delta>0$,
such that the following holds.
If at time $t\ge 0$, $x_{1,j}(t)=\nu_{j}$ for all $j\in \cj$,
and $x_{2,\ell}(t)\ge \epsilon$ for some fixed $\ell$,
then
$$
x_{1,\ell}(\tau) \ge \nu_{\ell}+\delta.
$$
\end{lem}

{\em Proof.}
(i) This easily follows from the definition of an FSP.
Clearly, shifted realizations of the primitive driving processes, 
defining FSP $x(\cdot)$, define $\theta_\tau x(\cdot)$.

(ii) If $x^n(\cdot)$ is a sequence of pre-limit trajectories defining FSP $x(\cdot)$,
then in a fixed small neighborhood of $t$, condition $\sum_j y^n_{0,j}(s) >0$ 
holds for all sufficiently large $n$. This means that (for large $n$), all new arrivals
in that neighborhood are routed to idle servers. Given the FSLLN properties
of driving trajectories, we easily obtain \eqn{eq-deriv-1}-\eqn{eq-deriv-2}
for any {\em regular} $t>0$. But then, given the continuity of $x(\cdot)$
and the fact that almost all time point are regular, we see that 
\eqn{eq-deriv-1}-\eqn{eq-deriv-2} must in fact hold for any $t$ 
(as long as $\sum_j y_{0,j}(s) >0$).

(iii) From (ii) we in particular have the following. For an FSP $x(\cdot)$,
at any $t\ge 0$ such that $\sum_j y_{0,j}(t) >0$ and $x_{2,j}(t)=0$ 
(i.e. $y_{1,j}(t)=x_{1,j}(t)$)
for all $j$,
$$
(d/dt) x_{k,j}(t) = 0, ~~k\ge 2, ~\forall j,
$$
$$
(d/dt) x_{1,j}(t) 
= \lambda (\beta_j - x_{1,j}(t)) /(\sum_{\ell}  (\beta_\ell - x_{1,\ell}(t))) 
- \mu_j  x_{1,j}(t).
$$
So, vector $(x_{1,j}(t), ~j\in \cj)=(y_{1,j}(t), ~j\in \cj)$ 
follows an ODE, which has unique 
solution, up to a point in time when 
$\sum_{j}  (\beta_j - x_{1,j}(t))=\sum_{j}  y_{0,j}(t)$ 
hits $0$.

(iv) By (ii) and the definition of $x^*$, $(d/dt)x(t)=0$ if $x(t)=x^*$.
Then we apply (iii).

(v) The FSP $x(\cdot)$, starting from the idle initial condition is unique 
 up to the first time $\tau_1$, at which 
$x_{1,j}(t)$ for one of the $j$ hits $\nu_j$.
From the structure of the ODE we observe that if $x_{1,j}(\tau_1)=\nu_j$
for one $j$, it has to hold for all $j$. Therefore, if $\tau_1<\infty$,
then $x(\tau_1)=x^*$. If so, by (i) and (iv), $x(t)=x^*$ for all $t\ge \tau_1$.
Then, by (iii), such FSP is unique; moreover,
\beql{eq-fsp-bound}
x(t) \le x^*, ~~t\ge 0.
\end{equation}

Consider now the sequence of {\em processes} $x^n(\cdot)$, 
starting from the idle initial state 
for each $n$. Uniqueness of the FSP starting from the idle initial condition,
along with Lemma~\ref{lem-ms-fsp-conv}, implies that
$x^n(\cdot)$ converges
(on the probability space constructed above in this section)
to this unique FSP: $x^n(\cdot) \to x(\cdot)$, u.o.c, w.p.1.
Recall that, for each $n$, process $x^n(\cdot)$ is stochastically monotone
non-decreasing (see \eqn{eq-zero-monotone}).
We conclude that the FSP $x(t), ~t\ge 0$, is non-decreasing in $t$.
Therefore, as $t\to\infty$, $x(t) \to x^{**}$ for some $x^{**} \le x^*$
(recall \eqn{eq-fsp-bound}). Finally, again from the structure of the ODE,
we see that $x^{**} = x^*$ must hold, because otherwise
$$
[(d/dt) \sum_j x_{1,j}(t)]_{x(t)=x^{**}} > 0.
$$
(vi) From (ii) and definition of $\nu_j$, using relation $y_{1,j}(t)=x_{1,j}(t)-x_{2,j}(t)$,
we have
$$
(d/dt) x_{1,j}(t) = \mu_j x_{2,j}(t), ~~j\in \cj.
$$
(For $t=0$ it is the right derivative.)
Also from (ii), we observe that 
in a sufficiently small fixed neighborhood of time $t$,
the expression for the derivative $(d/ds) x_{1,\ell}(s)$
must be uniformly Lipschitz continuous.
This implies that, for an arbitrarily small $\epsilon_1>0$,
in a (further reduced) small neighborhood $t$,
$(d/ds) x_{1,\ell}(s) \ge \mu_j \epsilon -\epsilon_1$;
which in turn implies the desired property.
$\Box$

\section{Proof of Theorem~\ref{thm-infinite}}
\label{sec-proof-infinite}

Since space $\cx$ is compact, any subsequence of $n$ has a further subsequence,
along which
\beql{eq-stationary-conv}
x^n(\infty) \Rightarrow x^{\circ}(\infty),
\end{equation}
where $x^{\circ}(\infty)$ is a random element in $\cx$.
Therefore, to prove Theorem~\ref{thm-infinite} it suffices to show that
any limit in \eqn{eq-stationary-conv} is equal (w.p.1) to $x^*$.

\begin{lem}
\label{lem-lower-bound}
Any subsequential limit $x^{\circ}(\infty)$ in \eqn {eq-stationary-conv}
is such that
$$
x^* \le x^{\circ}(\infty), ~~w.p.1.
$$
\end{lem}

{\em Proof.} For each $n$, consider the process $x^n(\cdot)$, 
starting from idle initial state. Consider any fixed $j$.
Fix arbitrary $\epsilon>0$, and choose $T>0$ large enough so that
the FSP $x(\cdot)$ starting from idle initial condition (as in
Lemma~\ref{lem-ms-fsp-prop}(v)) is such that $x_{1,j}(T) \ge \nu_j -\epsilon/2$.
Then, by Lemma~\ref{lem-ms-fsp-conv}, $\pr\{x_{1,j}^n(T) > \nu_j -\epsilon\} \to 1$.
We obtain 
$$
\liminf_{n\to\infty} \pr\{x_{1,j}^n(\infty) > \nu_j -\epsilon\} \ge 
\liminf_{n\to\infty} \pr\{x_{1,j}^n(T) > \nu_j -\epsilon\} = 1.
$$
Therefore, since $\{x_{1,j} > \nu_j -\epsilon\}$ is an open set, by the assumed 
convergence in distribution,
$$
\pr\{x_{1,j}^{\circ}(\infty) > \nu_j -\epsilon\} \ge 1.
$$
This holds for any $\epsilon>0$, so we have 
$\pr\{x_{1,j}^{\circ}(\infty) \ge \nu_j\} = 1$.
$\Box$

{\em Proof of Theorem~\ref{thm-infinite}.} 
First, we prove ergodicity (stability).
Let $x^n(\infty)$ be a random element, whose distribution is the lower 
invariant measure for the process $x^n(\cdot)$.
Consider the process, starting from the idle initial state, 
$x^n_{1,j}(0)=0, ~j\in\cj$.
Since $x^n(t)$ is stochastically monotone non-decreasing
and converges in distribution to $x^n(\infty)$ as $n\to\infty$, we observe that
the limit of the average expected (scaled) number of customer
service completions in $[0,T]$, as $T\to\infty$, is
$$
\lim_{T\to\infty} (1/T) \int_0^T [\E \sum_j \mu_j x_{1,j}^n(t)]dt 
= \E \sum_j \mu_j x_{1,j}^n(\infty).
$$
This limit cannot exceed $\lambda$, which is the the average expected 
(scaled) number of customer arrivals. (If the system initially has no customers,
the number of service completions in $[0,T]$ cannot, of course, exceed
the number of arrivals.) Therefore, 
\beql{eq-completion-rate}
\E \sum_j \mu_j x_{1,j}^n(\infty) \le \lambda.
\end{equation}
By Lemma~\ref{lem-unstable-measure}, for any $n$,
instability of the process is equivalent to condition 
\eqn{eq-instability-cond}, i.e.
$$
\pr \{x^n_{\infty,j}(\infty) = \beta_j\} =1,
$$
holding for at least one $j$.
Consider a subsequence of those $n$, for which the system is unstable,
with the above property holding for the same $j$. 
Consider a further subsequence, along which the convergence \eqn{eq-stationary-conv} to some
$x^{\circ}(\infty)$ holds; then, w.p.1, $x_{1,j}^\circ(\infty)=x_{\infty,j}^\circ(\infty) = \beta_j$
and (by Lemma~\ref{lem-lower-bound}) $x_{1,\ell}^\circ(\infty) \ge \nu_{\ell}$ for all $\ell$.
Therefore, along the chosen subsequence,
$$
\lim_n \E \sum_\ell \mu_\ell x_{1,\ell}^n(\infty) =
\E \sum_\ell \mu_\ell x_{1,\ell}^\circ(\infty)
\ge 
\beta_j \mu_j + \sum_{\ell \ne j} \nu_\ell \mu_\ell > \lambda.
$$
The contradiction with \eqn{eq-completion-rate} completes the proof of stability.

So, for every sufficiently large $n$, the process $x^n(\cdot)$ is stable,
and the lower invariant measure (which, by definition, is the distribution of $x^n(\infty)$) 
is its unique stationary distribution. 
Consider any subsequential limit $x^{\circ}(\infty)$  in \eqn{eq-stationary-conv},
long a subsequence of $n$; for  the rest of the proof, we consider $n$ along this subsequence.
By Lemma~\ref{lem-lower-bound},
$$
\E \sum_j \mu_j x_{1,j}^{\circ}(\infty) \ge \lambda.
$$
On the other hand, using \eqn{eq-completion-rate}, 
$$
\E \sum_j \mu_j x_{1,j}^{\circ}(\infty) = \lim_{n\to\infty} \E \sum_j \mu_j x_{1,j}^n(\infty) \le \lambda,
$$
and, therefore,
$$
\E \sum_j \mu_j x_{1,j}^{\circ}(\infty) = \lambda,
$$
which (again, recalling Lemma~\ref{lem-lower-bound}) is only possible 
when
\beql{eq-x1-is-nu}
x_{1,j}^{\circ}(\infty) = \nu_j, ~j \in \cj, ~~~w.p.1.
\end{equation}
It remains to show that 
\beql{eq-x2-is-0}
x_{2,j}^{\circ}(\infty) = 0, ~j \in \cj, ~~~w.p.1.
\end{equation}
Suppose not, that is for at least one $\ell$, 
$\pr\{x_{2,\ell}^{\circ}(\infty) > \epsilon\} = 2\epsilon_1$, for some 
$\epsilon>0$, $\epsilon_1>0$.
Then, for all sufficiently large $n$ (along the subsequence we consider), 
$\pr\{x_{2,\ell}^n(\infty) > \epsilon\} > \epsilon_1$. 
For each sufficiently large $n$, consider $x^n(\cdot)$ in steady-state,
that is $x^n(t)\stackrel{d}{=}x^n(\infty)$ for all $t\ge 0$.
Then $\pr\{x_{2,\ell}^n(0) > \epsilon\} > \epsilon_1$. Now, employing 
Lemma~\ref{lem-ms-fsp-conv} and Lemma~\ref{lem-ms-fsp-prop}(vi), we can easily
show that, for some $\tau>0$, $\delta>0$, and all large $n$,
$$
\pr\{x_{1,\ell}^n(\tau) \ge \nu_\ell + \delta/2\} > \epsilon_1/2.
$$
But then 
$$
\pr\{x_{1,\ell}^{\circ}(\infty) \ge \nu_\ell+ \delta/2\} \ge \limsup_{n\to\infty}
\pr\{x_{1,\ell}^n(\tau) \ge \nu_\ell + \delta/2\} \ge \epsilon_1/2,
$$
a contradiction with \eqn{eq-x1-is-nu}, which proves \eqn{eq-x2-is-0}.
$\Box$

\section{Generalizations}
\label{sec-generalization}

Our analysis relies mostly 
on the monotonicity property.
Monotonicity guarantees existence of the unique lower invariant measure
(for each scaling parameter $n$) for the process considered on the compactified
state space (whether or not the original process stochastically stable).
Then, proving stochastic stability 
and asymptotic optimality is essentially reduced to 
establishing the corresponding properties of 
the lower invariant measures.

Monotonicity property is preserved under various generalizations of our model.
We describe two of them in this section. In both cases, all our results
and proofs hold essentially as is.

\subsection{Queue-size dependent service rate}
\label{sec-parallel-service}

In our basic model we assumed that each server has a fixed processing
rate, independent of the queue length. This assumption is not realistic
 in many cases of interest. For example, a server may be
a processing ``device'' (physical or virtual) consisting
in fact of $C\ge 1$ independent ``sub-servers,'' that can work in parallel.
In this case, if the service rate of each sub-server is $\mu^{1}>0$,
the maximum processing rate $\mu=C\mu^1$ is achieved when there are
at least $C$ customers at the server, $Q\ge C$. The dependence
$f(Q)$ of the service rate on the queue length $Q$ is: $f(Q)=Q \mu^1$ when $Q<C$,
and $f(Q)=\mu$ when $Q\ge C$.

There may be other situations, where simultaneous service of multiple 
customers by a server is possible, but the services are not independent
(say, processing of different customers requires access to some shared
resources). In this case, the total service rate $f(Q)$ may be an
increasing function of $Q$, but increasing {\em sub-linearly}.

We now describe the model and PULL algorithm generalization, 
which accommodates the above considerations, while keeping 
the underlying Markov process a countable-state Markov chain,
and preserving monotonicity.
All results of this paper are easily extended to this generalized model.

The model is as before, except each server in pool $j$ has 
a more general service rate. For each $j$, there is a finite integer 
number $C_j$, $1 \le C_j \le B_j$, which is the server capacity,
in the sense of the maximum number of customers it can serve simultaneously.
The total service rate $f_j(Q)$, as a function of queue length $Q$,
is non-negative non-decreasing and such that $f_j(0)=0$
and $f(Q)=\mu_j$ for $Q\ge C_j$. We assume that the service requirement
of each customer is an independent exponentially distributed 
random variable with mean $1$. (This is consistent with 
the basic model considered in the paper.) The service
discipline in each server is arbitrary, as long as it is work-conserving and
non-idling.

The routing algorithm is generalized as follows. 
\begin{definition}[PULL algorithm generalization]
\label{def-pull-general}
At any time, if a server $i$ in pool $j$
has queue length $Q_i$, then the router has $\max\{C_j-Q_i,0\}$ 
pull-messages from this server. In other words,
at any time router has as many pull-messages from a server
as the server has available ``slots'' for additional customers to serve.
(A practical implementation of this, assuming pull-messages 
are never lost, is as follows.
When the server is ``initialized'', it sends $C_j$ pull-messages at once.
After that, the server sends one new pull-message immediately after
any service completion that leaves its queue length strictly less
than $C_j$.)
If at a customer arrival
the router has
available pull-messages (recall, that there may be multiple
pull-messages from any server), then it chooses one of them
uniformly at random, sends the customer to the corresponding server,
and destroys the ``used'' pull-message. If there are no available
pull-messages at a customer arrival, the customer is routed
uniformly at random to one of the servers in the system.
\end{definition}

Note that, as before, the router need {\em not} know anything 
about the parameters or the current states of the servers, besides the current set of
available pull-messages. Again, from the router's point of view
all servers form a single pool, despite possible differences
in the servers' parameters.

The queue length process for this model and PULL algorithm
is a monotone countable-state-space Markov chain.
All our results and proofs easily generalize.

\subsection{More general service time distributions}
\label{sec-dhr}

The assumption that the service times have exponential distribution,
can also be relaxed. To simplify the discussion, let us assume for now
that, as in the basic model, each server is a ``single-server'' (has 
constant processing speed, regardless of the queue length),
employing FCFS discipline.

Assume that the service time distribution
in each pool $j$ has {\em decreasing hazard rate} (DHR), 
and has positive finite mean $1/\mu_j$.
A distribution 
on $\R_+$, with complementary distribution function $F^c(z), ~z\ge 0$,
has DHR if the hazard rate
$$
- \frac{(d/dz) F^c(z)}{F^c(z)}
$$
is a non-increasing function of $z$. Exponential distribution 
with mean $1/\mu$ is a special case,
with constant hazard rate $\mu$. Another important example is the 
(heavy-tailed)
Pareto distribution:
$$
F^c(z) = [1+\sigma z]^{-\alpha},
$$
with parameters $\sigma>0$ and $\alpha>1$; it has finite mean value
$\mu^{-1}=[\sigma(\alpha -1)]^{-1}$.
{\em If service time distributions have DHR, then the assumption that 
the service 
in each queue is FCFS order is essential.}
The state of queue $i$ is the pair
$(Q_i, H_i)$, where, as before, $Q_i\ge 0$ is the (integer) queue length
and $H_i \ge 0$ is the (real) elapsed service time of the head-of-the-line
customer. (If $Q_i=0$ then necessarily $H_i=0$.)
The order $(Q_i, H_i) \le (Q'_i, H'_i)$ is understood component-wise.

The compactification of the state space $\Z_+ \times \R_+$ 
of one server in pool $j$ is done in two steps.
In the first step, we compactify $\Z_+ \times \R_+$ to
$\bar \Z_+ \times \bar \R_+$, where each component
$\bar \Z_+ = \Z_+ \cup \{\infty_Q\}$ 
and $\bar \R_+ = \R_+ \cup \{\infty_H\}$ is compactified separately
(where $\infty_Q$ and $\infty_H$ are the corresponding ``points at infinity"),
with the product topology on $\bar \Z_+ \times \bar \R_+$.
The second step depends on whether the minimum
hazard rate
$$
\gamma_j \doteq \lim_{z\to\infty} \left[- \frac{(d/dz) F^c_j(z)}{F^c_j(z)}\right]
$$
is zero or not. (Here $F^c_j(\cdot)$ is the complementary 
distribution function of a service time in pool $j$.)
If $\gamma_j>0$, we further identify all points
$(Q_i, H_i)$ with $Q_i=\infty_Q$ as a single point $\infty$ at infinity;
if $\gamma_j=0$, we further identify all points
$(Q_i, H_i)$ with either $Q_i=\infty_Q$ or $H_i=\infty_H$ 
as a single point $\infty$ at infinity.
The server state $(Q_i, H_i)=\infty$ is such that it never changes --
neither service completions nor new arrival to the server affect it.
The order relation is naturally extended to the compactified state space.

The Markov process, describing system evolution, is monotone.
Its stability is understood more generally, as positive Harris recurrence,
and is equivalent to the fact that the lower invariant measure is proper, i.e.,
almost surely every server state belongs to $\Z_+ \times \R_+$ .

The corresponding mean field (fluid-scaled) processes
and fluid sample paths in this model are more general -- 
the state component for each $(k,j)$
 is not just a number, but a function describing the distribution 
of elapsed service times among the servers in pool $j$ with queue length $k$.
The equilibrium point is defined accordingly; its projection on space $\cx$,
describing queue lengths only (without regard to elapsed service times),
is still $x^*$ as defined in \eqn{eq-invar-point1}-\eqn{eq-invar-point2}
 -- it is invariant w.r.t. service time 
distributions given their means $1/\mu_j$. 
The appropriately generalized version of Theorem~\ref{thm-infinite} holds 
under these assumptions, with essentially same proof.

The model can be further generalized to assume that each 
server in pool $j$ consists of a finite number
$C_j \le B_j$ of ``sub-servers'' that
can work independently in parallel (as was described at the beginning
of Section~\ref{sec-parallel-service}). 
Within each server, the customers are 
allocated to sub-servers in FCFS order. (This is essential.)
The service time distribution of a customer in one sub-server
in pool $j$ has DHR with mean $C_j /\mu_j$; 
so that the maximum processing rate is $C_j [C_j /\mu_j]^{-1}=\mu_j$.
The PULL algorithm is as in Definition~\ref{def-pull-general}.
The state of a server, besides the queue length, will now contain 
the elapsed service times of the customers in service;
the states equal up to a permutation of sub-servers are identified;
the state space is compactified analogously to the way it is done above
for the single-server case; the natural order relation is considered.
The corresponding Markov process is monotone.
Theorem~\ref{thm-infinite}
 generalizes to this model as well and, again,
it implies that asymptotically, under the subcritical load condition
\eqn{eq-load}, the steady-state probabilities of waiting or blocking, vanish.

\bibliographystyle{acmtrans-ims}
\bibliography{biblio-stolyar}


\appendix

\section{Additional corollaries from the main results}
\label{app1}

\subsection{System with infinite buffers: insensitivity to queueing disciplines at the servers.}

It is described in Section~\ref{sec-dhr} how our main result,  Theorem~\ref{thm-infinite}, is generalized to service time distributions with DHR,
under the FCFS assumption on the queueing discipline at each server. The FCFS assumption is essential for our approach to work. However, in the special case when 
{\em all buffer sizes are infinite and service rates are independent of the queue length  (i.e. 
$B_j=\infty$ and $C_j=1$ for all $j$),} the FCFS assumption is, in fact, {\em not} essential. 
In this special case, no arriving customer is ever blocked and the unfinished work at each server ``drains'' at unit rate (when it is non-zero) under any service discipline that
is work-conserving and non-idling. By the definition of PULL,
 the assignment of each arriving customer depends only on which servers are currently idle, i.e.,  which of them have zero unfinished work. Therefore, we obtain the
following simple

\begin{lem}
\label{lem-insensitivity}
Consider the system with infinite buffer sizes and queue length independent service rates. The customer service time distributions in different pools are arbitrary.
Define a server state as its total amount of unfinished work, and the system state accordingly. Then, under PULL algorithm, the system state process is 
invariant with respect to the service discipline at each server, as long as it is non-idling and work-conserving.
\end{lem}

As a corollary of the argument in Section~\ref{sec-dhr} (which is for the FCFS discipline at the servers) and Lemma~\ref{lem-insensitivity}, we see that {\em in the special case of  infinite buffer sizes and queue length independent service rates,
the extension of Theorem~\ref{thm-infinite} to DHR is valid for arbitrary non-idling work-conserving disciplines at the servers.}

\subsection{Service time distributions with a positive lower bound on the hazard rate.}

Suppose, for each $j$ the service time distribution is such that its hazard rate is lower bounded by $\gamma_j>0$. Then we can use the monotonicity approach to compare this system
to the corresponding system with exponential service time distributions with rates $\gamma_j$. Namely, using the same constructions and arguments as in 
 Section~\ref{sec-dhr}, we can easily verify the following 

\begin{lem}
\label{lem-hazard-lower}
Consider the system described in Section~\ref{sec-dhr},
with FCFS service discipline at each server. 
The service time distribution in pool $j$ has the hazard rate lower bounded by $\gamma_j>0$.
Let us label this system by S1. Consider the corresponding system with exponential service time distributions with rates $\gamma_j$; 
let us label this system by S2. Then, if initial state of S1 is dominated by that of S2, the processes for the two systems can be coupled so that, w.p.1,
this dominance relation prevails at all times.
\end{lem}

As a corollary of  Theorem~\ref{thm-infinite} and Lemma~\ref{lem-hazard-lower}, we obtain the following

\begin{prop}
\label{prop777}
Suppose condition 
\beql{eq-stabil-gamma}
\lambda < \sum_j \beta_j \gamma_j
\end{equation}
holds (which is condition \eqn{eq-load} for system S2). Then,
for all sufficiently large $n$, the system S1 state process is positive Harris recurrent and its unique stationary distribution is stochastically dominated by that of
system S2. 
\end{prop}

Proposition~\ref{prop777} implies that the asymptotic optimality of PULL prevails for system S1 under condition \eqn{eq-stabil-gamma}
and FCFS discipline at each server.

\end{document}